# Some improvements on the one-step inverse isogeometric analysis by proposing a multi-step inverse isogeometric methodology in sheet metal stamping processes

Amir Reza Isazadeh, Mansoor Shamloofard, Ahmad Assempour[*]

**Abstract**

The isogeometric methodology has been successfully implemented in one-step inverse analysis of sheet metal stamping processes. However, these models are not capable of analyzing forming processes that require severe deformation and several forming stages. This paper presents a multi-step inverse isogeometric methodology to enhance the precision of one-step models in predictions of the initial blank, strain distributions, and drawability of the formed parts. This methodology deals with the minimization of potential energy, deformation theory of plasticity, and considering membrane elements. The presented methodology utilizes the NURBS basis functions to create the final, middle, and blank geometries and also to analyze sheet metal deformation. The characteristics of the applied formulations make it possible to simultaneously observe the effects of changing part parameters on its formability. One advantage of this approach is that the linear system of governing equations is solved without concerning about the convergence. Besides, the presented methodology can successfully generate the middle geometry and to restrict the movements of physical nodes along the middle surface, by presenting a new NURBS-based mapping and sliding constraint technique. The performance of the presented model is experimentally and numerically evaluated under two classical problems, including the forming of a rectangular box and a two-step drawing of a circular cup. Results comparisons indicate the credibility of the presented model in prediction of forming parameters at a low computation time.

[*] Center of Excellence in Design, Robotics and Automation, Department of Mechanical Engineering, Sharif University of Technology, Azadi Avenue, P.O. Box 11365-9567, Tehran, Iran
e-mail: assem@sharif.edu

## 1. Introduction

Numerical simulation of sheet metal stamping processes is a nonlinear problem that involves nonlinearity of geometry, material, and contacts. Therefore, incremental methods, which use plastic flow theory, require tremendous computation time and powerful computer resources. As a result of these properties, incremental-based models are extensively employed in the final stages of design. On the other hand, inverse approaches utilize the deformation theory of plasticity and proportional loading conditions, which result in less computational cost and make this methodology well-suited for the initial stages of design.

Early research studies to develop one-step inverse approaches were conducted by Majlessi and Lee [1-2]. Their nonlinear model was based on minimum potential energy and could not be solved without considering blank holder and friction forces. Batoz and Guo [3-4] developed an inverse finite element method (IFEM) based on the principle of virtual work. This inverse model could be implemented in the simulation of stamping processes without blank holder and friction forces, in contrast to the inverse models presented in [1-2]. For a large number of applications, membrane effects are dominant; however, bending effects should be included into the computation to enhance the accuracy of the forming analysis, especially in the case of increasing sheet thickness to punch radius ratio. Batoz et al. [5] presented a rotation-free shell element in the content of the inverse approach for considering the bending effects. As the mentioned inverse models [1-5] lead to nonlinear governing equations, the convergence of these methods strongly depends on the provided initial solution. To overcome this limitation, Liu and Karima [6], Assempour et al. [7] proposed an unfolding technique in which the 3D problem is transformed into a 2D space, and the governing equations are solved without any concern about the convergence.

Although one-step inverse models provide reasonable results with less computation time in many forming processes, these models cannot accurately analyze forming processes that require severe deformation and several forming stages. To enhance the accuracy of one-step models in these conditions, multi-step inverse approaches have been presented in several research studies. For instance, Majlessi and Lee [8] developed a multi-step inverse finite element model that was applicable only for axisymmetric problems. Then, Lee and Cao [9] added bending effects to this model by considering the rotation-free shell element. Lee and Huh [10] introduced an inverse mapping to estimate the initial guess of the intermediate configuration. This method is not capable of analyzing complex geometries. To enhance the credibility of this model, Kim and Huh [11-12] proposed a new direct mesh mapping. Huang et al. [13] presented a general method for non-axisymmetric geometries, called the modified arc-length search method. They applied the concept of sliding constraint surface to restrict the displacement of nodes along the intermediate surface. Furthermore, Bostanshirin et al. [14] introduced a new mapping technique to analyze multi-stage forming processes.

The traditional one-step and multi-step inverse models have been formulated based on the finite element method (FEM). Despite numerous advantages of this method, a large percentage of the required time for FEM analysis is devoted to the mesh generation stage [15]. In addition, in analysis of a complicated part, this time-consuming process is increased due to the separate nature of modeling and analysis representations in this approach.

To eliminate spending time on mesh generation, isogeometric analysis (IGA) was proposed by Hughes et al. [16]. This methodology takes advantage of Non-Uniform Rational B-Spline (NURBS) basis functions to draw curves and surfaces accurately as well as to approximate field variables. In IGA, the coupling of computer-aided design (CAD) and computer-aided engineering (CAE) eliminates the mesh generation time. This tight integration between modeling and analysis representations makes this methodology considerably beneficial for shape optimization problems [17-21], and due to higher-order continuity of NURBS basis functions, many research studies have utilized this method for analysis of shell and plate problems [22-27]. Contact modeling [28-29], simulation of fluid-solid interaction [30-32], and structural vibration analysis [33-36] are other fields of research that IGA has successfully implemented. Compared to FEM, one drawback in the NURBS-based IGA is the inability to localize the mesh refinement. To solve this problem, many alternative spline-based techniques such as T-splines [37] and hierarchical B-splines [38] have been developed.

Forward incremental analysis of sheet metal forming in the isogeometric framework has been recently studied by Benson et al. [39]. They developed a NURBS-based isogeometric shell formulation and used it for simulation of sheet metal stamping processes. Similar to this research, Ambati et al. [40] presented a stress-based isogeometric thin shell formulation to consider large deformations.

More recently, the inverse isogeometric formulation has been developed to study sheet metal forming simulations. For instance, a one-step inverse approach in the framework of IGA was developed by Zhang et al.[41-42]. Using the deformation theory of plasticity and considering the principle of virtual work are the foundations of this inverse IGA approach. Their nonlinear governing equations are solved using the Newton-Raphson iterative algorithm. To ensure the convergence of this method, they used an energy-based algorithm [43] to provide an appropriate initial solution. Furthermore, Shamloofard et al.[44] presented an inverse isogeometric formulation by employing a new approach to transfer the 3D problem to a two-dimensional space. This approach ensures convergence and does not require an initial solution. They also proposed a new non-uniform friction model.

According to the author's knowledge, multi-step inverse isogeometric analysis (IIGA) has not been studied so far. The main objective of this research study is to develop a new multi-step IIGA to analyze sheet metal stamping processes requiring several forming stages and/or severe plastic deformation. Using the isogeometric approach, both geometry modeling and forming analysis are simultaneously carried out using the NURBS basis functions. In the proposed methodology, a new NURBS-based mapping is presented to obtain the first estimation for representation of the intermediate geometry. In addition, a method based on the concept of sliding constraint is utilized to control the movements of the nodes along the intermediate surface. Similar to the assumptions detailed in Ref [44], the deformation theory of plasticity and minimization of potential energy lead to linear governing equations in the presented model. In order to evaluate the accuracy of the proposed method, sheet metal forming analysis of a rectangular box with high drawing depth and also two-stage drawing of a circular cup are studied, and the results of the present model have been verified experimentally and numerically.

## 2. Preliminary study
### 2.1. Geometric modeling

B-spline functions are the most common parametric functions for modeling of curves and surfaces. In B-splines, each coordinate is defined as a function of independent variables. Non-uniform rational B-spline (NURBS) is the generalized form of B-splines which considers weights for control points. A knot vector and a series of control points are needed to define a curve by NURBS basis functions. A knot vector consisting of real numbers in the parametric space is expressed as $\Xi = \{\xi_1, \xi_2, \dots, \xi_{i+p+1}\}$ where $\xi_i \in R$ is the $i_{th}$ knot and $i = 1, 2, \dots, n+p+1$ is the knot index. Also, $p$ and $n$ represent the polynomial degree of basis functions and the number of control points.

In a single dimension, NURBS basis functions are defined by piecewise rational functions as follows:

$$R_i^p(\xi) = \frac{N_{i,p}(\xi)w_i}{W(\xi)} = \frac{N_{i,p}(\xi)w_i}{\sum_{i=1}^{n} N_{i,p}(\xi)w_i} \tag{1}$$

Where, $w_i$ and $N_{i,p}$ represent the weight of $i_{th}$ control point and B-spline basis functions which are calculated using the Cox–de Boor formula [45]:

$$N_{i,0}(\xi) = \begin{cases} 1 & \text{if } \xi_i \leq \xi < \xi_{i+1} \\ 0 & \text{otherwise} \end{cases}$$

$$N_{i,p}(\xi) = \frac{\xi - \xi_i}{\xi_{i+p} - \xi_i} N_{i,p-1}(\xi) + \frac{\xi_{i+p+1} - \xi}{\xi_{i+p+1} - \xi_{i+1}} N_{i+1,p-1}(\xi) \tag{2}$$

Considering NURBS basis functions and corresponding control points, a NURBS curve is expressed as follows [46]:

$$C(\xi) = \sum_{i=1}^{n} R_i^p(\xi) B_i \tag{3}$$

In which, $B_i$ are control points. Also, using a control net $\{B_{i,j}\}$, $i = 1,2,...,n, j = 1,2,...,m$ and two knot vectors ($\Xi = \{\xi_1, \xi_2, ..., \xi_{n+p+1}\}$ and $\mathcal{H} = \{\eta_1, \eta_2, ..., \eta_{m+q+1}\}$), a NURBS surface tensor product is constructed as below:

$$S(\xi, \eta) = \sum_{i=1}^{n} \sum_{j=1}^{m} R_{i,j}^{p,q}(\xi, \eta) B_{i,j} \tag{4}$$

$$R_{i,j}^{p,q}(\xi, \eta) = \frac{N_{i,p}(\xi) M_{j,q}(\eta) w_{i,j}}{\sum_{i=1}^{n} \sum_{j=1}^{m} N_{i,p}(\xi) M_{j,q}(\eta) w_{i,j}^n} \tag{5}$$

Where, $N_{i,p}(\xi)$ and $M_{j,q}(\eta)$ are the B-spline basis functions corresponding to $\mathcal{H}$ and $\Xi$ knot vectors.

### 2.2. Isogeometric analysis

Due to the valuable properties of the NURBS basis functions (described in Ref. [16]), the NURBS-based IGA takes advantage of these functions to draw geometry as well as to approximate field variables. The isoparametric concept is employed in this approach, and an exact geometric map is implemented, in contrast to the approximate geometry map used in the isoparametric FEM. Using this definition, a field variable ($u$) is approximated as the following equation:

$$u(\xi, \eta) \approx \sum_{i=1}^{n} R_i^{p,q}(\xi, \eta) u_i \tag{6}$$

In which, $u_i$ represent the control variables and $R_i^{p,q}(\xi, \eta)$ are the NURBS basis functions.

### 2.3. One-step transfer-based IIGA [44]

In the design of sheet metal stamping process using inverse analysis, the geometry of the final configuration ($C^F$) and the thickness of the initial sheet ($C^I$) shown in Fig. 1 are the inputs. The main idea is to predict the shape and size of the initial blank (positions of control points in the initial flat sheet) as well as to calculate thicknesses, strains and stresses of the final configuration ($C^F$), according to Fig. 1. For this purpose, the transfer-based one-step IIGA employs the deformation theory of plasticity and the principle of minimum potential energy. In this methodology, the following items are considered to develop the governing equations [44]:

- ✓ Only the initial and final configurations are considered;
- ✓ Interaction forces and contact modeling between blank-die and blank-punch are not included into the computations. In this method, in-plane forces are considered as forces which transform the initial configuration to the final part;
- ✓ Hill yield criterion under the plane stress assumption is used to update the equations;
- ✓ The planar anisotropic condition using Lankford coefficients is employed to define anisotropic properties.

In the transfer-based IIGA, a 3D problem is primarily transformed into a 2D space using the isogeometric transfer concept [44]. In this iterative methodology, the governing equations are solved in two-dimensional space, without concerning about the pre-estimation solution, and material properties are updated in each iteration step. Therefore, less computation time is required in this method, and the convergence of the solution is guaranteed. Details of the IG transfer are given in Ref. [44].

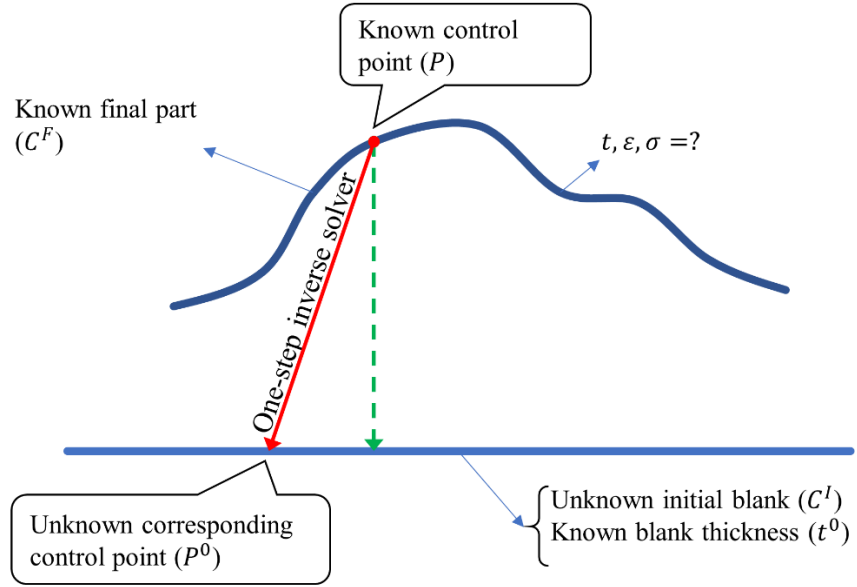

Fig. 1. Description of the one-step inverse method

## 3. Multistep transfer-based IIGA

One-step transfer-based IIGA formulation proves its ability in the analysis of sheet metal stamping processes. However, a multistep transfer-based approach should be offered in some cases in which several stages of forming and also severe plastic deformation are required. In this section, a multi-step transfer-based IIGA methodology is proposed to improve the accuracy of the one-step transfer-based method. This new multi-step IIGA model follows the basic assumptions of the transfer-based one-step IIGA approach and considers initial, final and intermediate stages of the forming process.

The geometries of the final part, intermediate punch and die, thickness and material properties of the initial sheet, friction coefficient, and the blank holder force are the inputs of this approach. Also, its outputs include the shape and size of the initial blank, geometry of the middle configuration, and thicknesses/strains/stresses of the middle and final parts. In this method, a configuration is initially estimated as the middle part, using a NURBS-based mapping. Then, by employing the principle of minimum potential energy, the initial solution is updated, and the strains and stresses are calculated. Next, the one-step transfer-based IIGA is applied to the updated middle configuration, and the shape of the blank, strains, and stresses of the middle configuration are obtained. Eventually, the material properties of the final part are modified, and this process is repeated up to the convergence.

## 3.1. Governing equations

Similar to the one-step transfer-based IIGA, the presented approach is based on iterations. In the first step of this process, the first estimation of the middle configuration is constructed to evaluate the middle configuration. Then, in the next iterations, this configuration is updated by applying the principle of minimum potential energy on the final configuration and using isogeometric membrane elements which lead to the following equation [44]:

$$\sum_{i=1}^{ne} t_i \int_{A_i^T} [B_i]^T [D_i][B_i]\{d_i\}dxdy - \sum_{i=1}^{ne} F_i = 0$$
$$\sum_{i=1}^{ne} t_i \int_{A_i^T} [B_i]^T [D_i][B_i]\det(J_1)_i d\xi d\eta\{d_i\} - \sum_{i=1}^{ne} F_i = 0 \quad (7)$$
$$\sum_{i=1}^{ne} t_i \int_{-1}^{1} [B_i]^T [D_i][B_i]\det(J_1)_i \det(J_2)_i dsdt\{d_i\} - \sum_{i=1}^{ne} F_i = 0$$

In which, $t_i$, $[B]$, $[D]$, $\{d_i\}$, $ne$, $F_i$, $(x,y)$, $(\xi,\eta)$ and $(s,t)$ represent the thickness of each element, the membrane strain–displacement matrix, material property matrix, the displacement vector of the control points, number of isogeometric elements, external forces, physical space, parametric space, and integration space, respectively. Also $J_1$ and $J_2$ are related to the conversion of physical space to parametric space and conversion of parametric space to integration space, as shown in Fig. 2. NURBS basis functions are employed for the first conversion, and polynomial function is used for the second.

The strain–displacement matrix is defined as below:

$$[B] = [L][R] = \begin{bmatrix} \frac{\partial}{\partial x} & 0 \\ 0 & \frac{\partial}{\partial y} \\ \frac{\partial}{\partial y} & \frac{\partial}{\partial x} \end{bmatrix} \begin{bmatrix} R_1 & \dots & R_{ncp} & 0 & \dots & 0 \\ 0 & \dots & 0 & R_1 & \dots & R_{ncp} \end{bmatrix} \quad (8)$$

In which, $ncp$ is the number of control points in each isogeometric element.

In Eq. (7), $F_i$ represent external forces which are obtained by combinations of in-plane and friction forces as the following:

$$\vec{F} = \vec{F}_f + \vec{F}_{in-plane} \quad (9)$$

In which, $\vec{F}_f$ is the non-uniform friction force vector, which is calculated using the equations presented in Ref. [44]. Also, $\vec{F}_{in-plane}$ is the in-plane force vector which transforms the middle configuration into the final network. In order to calculate the in-plane forces, displacements of control points between these configurations in each element and also the stiffness matrices of final part elements are calculated as the following equations:

$$\{\Delta X\}_L^e = \{X\}_L^M - \{X\}_L^F \tag{10}$$

$$[k]_L^e = t_e \int_{-1}^{1} [B_e]^T [D_e][B_e] det(J_1)_e det(J_2)_e ds dt \tag{11}$$

Where, $\{X\}_L^M$ and $\{X\}_L^F$ are the local coordinates of control points in middle and final configurations, respectively. Considering local stiffness matrices ($[K]_L$), the in-plane force vector is computed in the local and global coordinates using the translation matrix $[T]$ which will be discussed in the next subsection.

$$\{\vec{f}_{in-plane}\}_L^e = [k]_L^e \{\Delta X\}_L^e \tag{12}$$

$$\{\vec{f}_{in-plane}\}_G^e = [T]^T \{\vec{f}_{in-plane}\}_L^e \tag{13}$$

Assembling the in-plane forces in the global coordinate system (CSYS) results in the total in-plane force vector as follows:

$$\vec{F}_{in-plane} = \sum_{e=1}^{n} \{\vec{f}_{in-plane}\}_G^e \tag{14}$$

To assemble the calculated stiffness matrices, local stiffness matrix of each element should be specified in the global CSYS as the following equation:

$$[k]_G^e = [T]_F^T [k]_L^e [T]_F \tag{15}$$

### 3.1.1. Local and global coordinate systems

In Fig.2, integration, parametric, and physical spaces of a NURBS surface created by two knot vectors {0,0,0,0.5,1,1,1} and {0,0,0,0.3,0.6,1,1,1} are displayed. As shown in this figure, knot vectors discretize the surface to NURBS elements. Considering parametric points of an element $((\xi_i, \eta_j), (\xi_{i+1}, \eta_j), (\xi_{i+1}, \eta_{j+1}), (\xi_i, \eta_{j+1}), (\frac{\xi_i + \xi_{i+1}}{2}, \frac{\eta_i + \eta_{i+1}}{2}))$, physical nodes of this element ($P_1$, $P_2$, $P_3$, $P_4$, $P_C$) can be calculated using Eq. (4). Furthermore, unit vectors of the local CSYS in each element can be defined as follows:

$$\vec{e}'_1 = \frac{\vec{P}_2 - \vec{P}_1}{\|\vec{P}_2 - \vec{P}_1\|}$$

$$\vec{e}'_3 = \frac{(\vec{P}_2 - \vec{P}_1) \times (\vec{P}_4 - \vec{P}_1)}{\|(\vec{P}_2 - \vec{P}_1) \times (\vec{P}_4 - \vec{P}_1)\|} \tag{16}$$

$$\vec{e}'_2 = \vec{e}'_3 \times \vec{e}'_1$$

Transformation matrix between two local and global coordinates can be expressed as:

$$Q_{ij} = cos(\vec{e}_i, \vec{e}'_j) \tag{17}$$

In which, $\vec{e}_i$ indicate the unit vectors of the global CSYS. Therefore, control points of each element can be transferred to their local CSYS as follows:

$$\{X_i\}^e_L = [Q](\{X_i\}^e_G - \{P_C\}^e_G) \tag{18}$$

In which, $\{X_i\}^e_G$, $\{X_i\}^e_L$ represent the coordinates of $i_{th}$ control point in the global and local CSYS, respectively.

Finally, the translation matrix between local and global CSYS is deduced as:

$$[T] = \begin{bmatrix} Q_{11}.I_{ncp} & Q_{12}.I_{ncp} \\ Q_{21}.I_{ncp} & Q_{22}.I_{ncp} \\ Q_{31}.I_{ncp} & Q_{32}.I_{ncp} \end{bmatrix} \tag{19}$$

In which, $I_{ncp}$ and $Q_{ij}$ indicate a $ncp$-by-$ncp$ identity matrix and the transformation matrix between the global and local CSYS in the middle configuration.

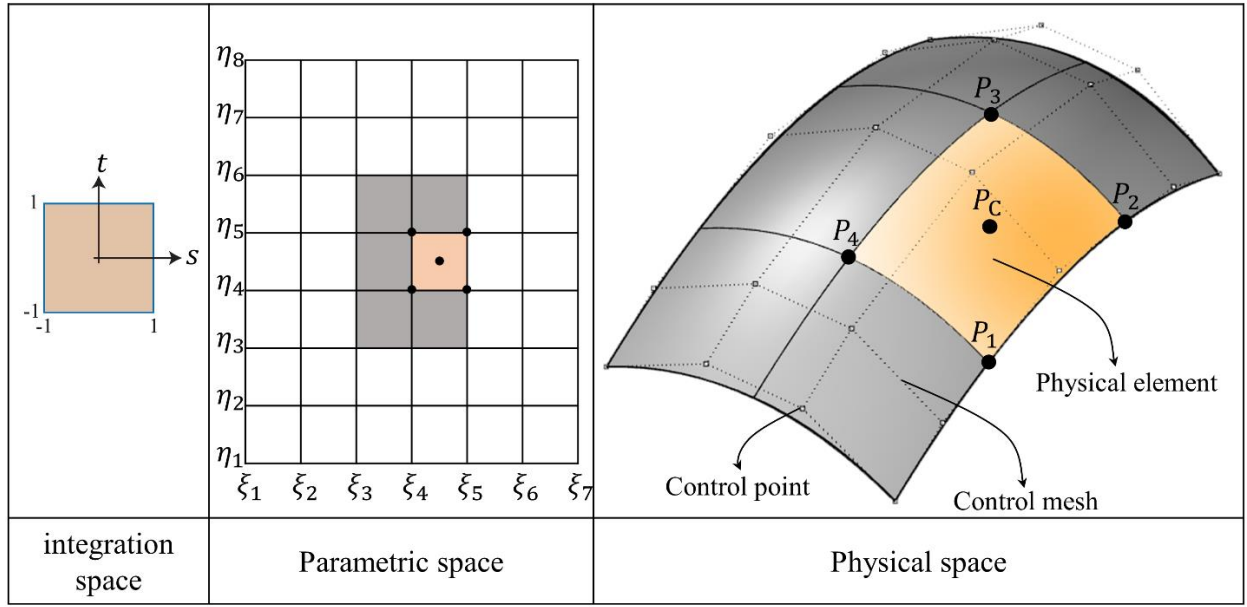

Fig. 2. Different spaces in IGA

### 3.2. Obtaining the middle configuration using NURBS-based mapping

In this research, the first estimation of the middle configuration ($C^M_{mapping}$) is generated using a new NURBS-based mapping. Prior to constructing this estimation, a user-defined middle configuration ($C^M_{user}$) is created considering the punch and die profiles of this stage, which are assumed to be known. The only constraint to construct the user-defined middle configuration is that the initial blank of the final configuration should be surrounded by that of the user-defined middle configuration. This user-defined configuration is only used to generate the first estimation of the middle part by the NURBS-based mapping. The estimated middle part in the first iteration ($C^M_{mapping}$) will be updated in the next iterations by satisfying the minimum potential energy described in section 3.1. Fig. 3 illustrates how the last version of the middle part is specified in the presented methodology. The procedure is categorized as the following five stages:

1- The one-step IIGA is initially applied to $C^F$ and $C^M_{user}$. As a result, the initial blanks of the final part ($B^F$) and user-defined middle surface ($B^M_{user}$) are obtained, according to Fig. 3(a).

2- As shown in Fig. 3 (a), a physical point (A) is located in $C^F$ and $B^F$. The isogeometric element of the $B^M_{user}$ incorporating this point is found using a searching algorithm.

3- Considering the point (A) as a physical point on the $B^M_{user}$, the parametric coordinates ($\xi, \eta$) are specified using the relationship between the physical point and control points of the incorporating element as follows:

$$\begin{Bmatrix} x(\xi,\eta) \\ y(\xi,\eta) \\ z(\xi,\eta) \end{Bmatrix} = \begin{bmatrix} R_1 & \cdots & R_{ncp} & 0 & \cdots & 0 & 0 & \cdots & 0 \\ 0 & \cdots & 0 & R_1 & \cdots & R_{ncp} & 0 & \cdots & 0 \\ 0 & \cdots & 0 & 0 & \cdots & 0 & R_1 & \cdots & R_{ncp} \end{bmatrix} \\ \times \{X_1, \dots, X_{ncp}, Y_1, \dots, Y_{ncp}, Z_1, \dots, Z_{ncp}\}^T \quad (20)$$

In which, $(x, y, z)$ and $(X, Y, Z)$ represent the coordinates of physical point and control points.

4- The corresponding physical point in the $C^M_{user}$ is specified using the calculated parametric coordinates, according to Fig. 3 (b).

5- This procedure is repeated for all physical nodes of $C^F$ and a new middle configuration is obtained (Fig. 3 (c)). Finally, by applying the governing equations presented in Section 3.1 up to the convergence, the middle part achieved by the mapping is updated, and the last version of the middle part will be specified.

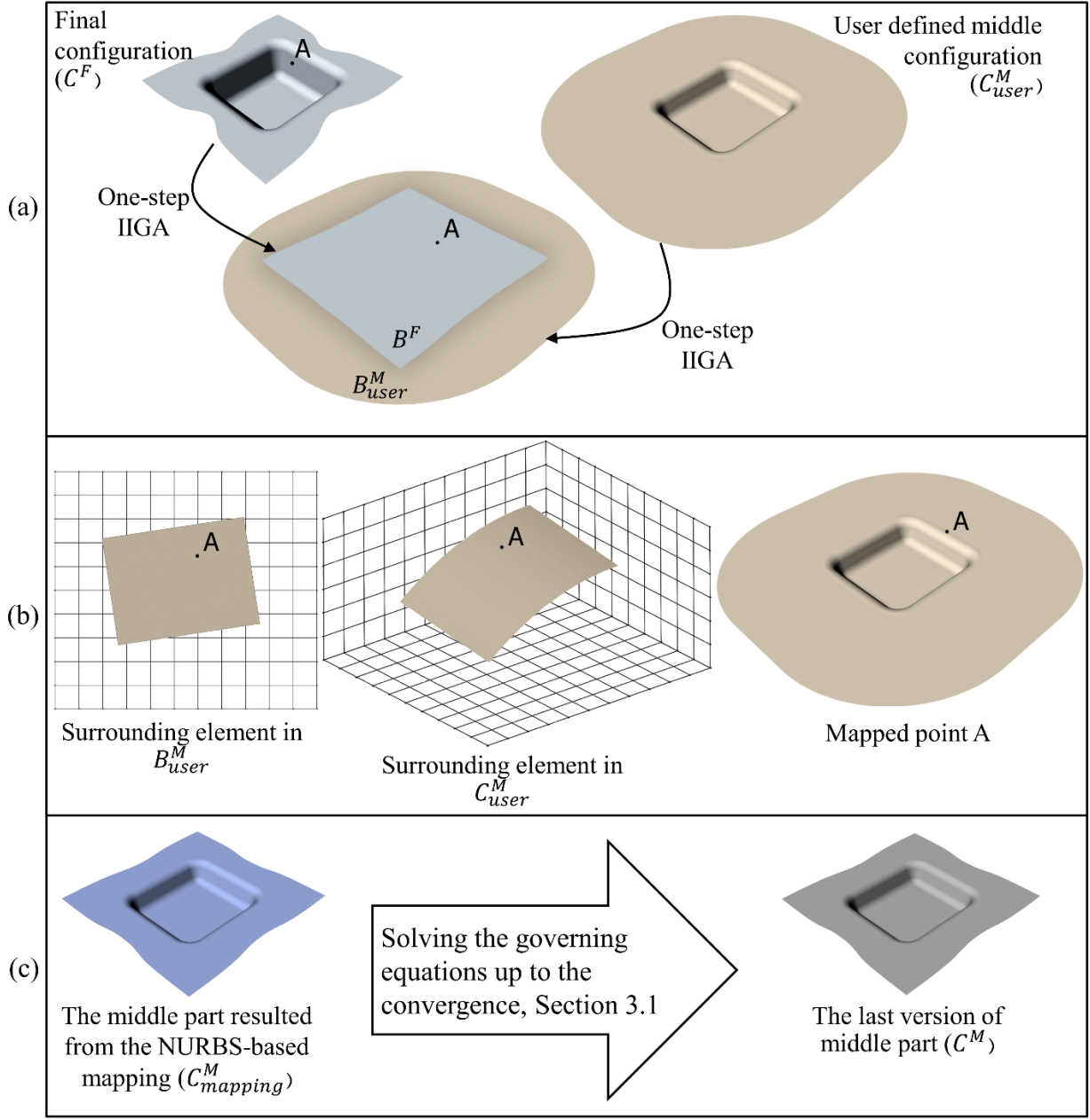

Fig. 3. Description of the different middle configurations used in the present multi-step method

### 3.3. Concept of sliding constraint surface

To restrict the movements of physical nodes along the surface of the intermediate part, a new constraint is used in this research. To do so, local coordinates that are perpendicular to the middle surface for each control point are initially defined, as illustrated in Fig. 4 (a). Then, the global stiffness matrix, the displacement vector, and the force vector are transformed from global CSYS into the new local CSYS. By eliminating the movements of control points along the normal direction of the surface (Fig. 4 (a)), the linear system of equations is solved. One shortcoming regarding this approach is related to applying large displacements to the middle configuration. In this case, the control points will move tangent to the middle surface, especially in the die and punch profile radii, as depicted in Fig. 4 (b). Therefore, the problem cannot converge. To address this limitation, the computed displacements are divided into several small displacements. The first part of the divided displacement is initially applied, and the normal coordinate systems are recalculated. In the subsequent stages, the remained parts of displacements will be separately considered up to the last part. Therefore, a continuous path is obtained at the middle surface.

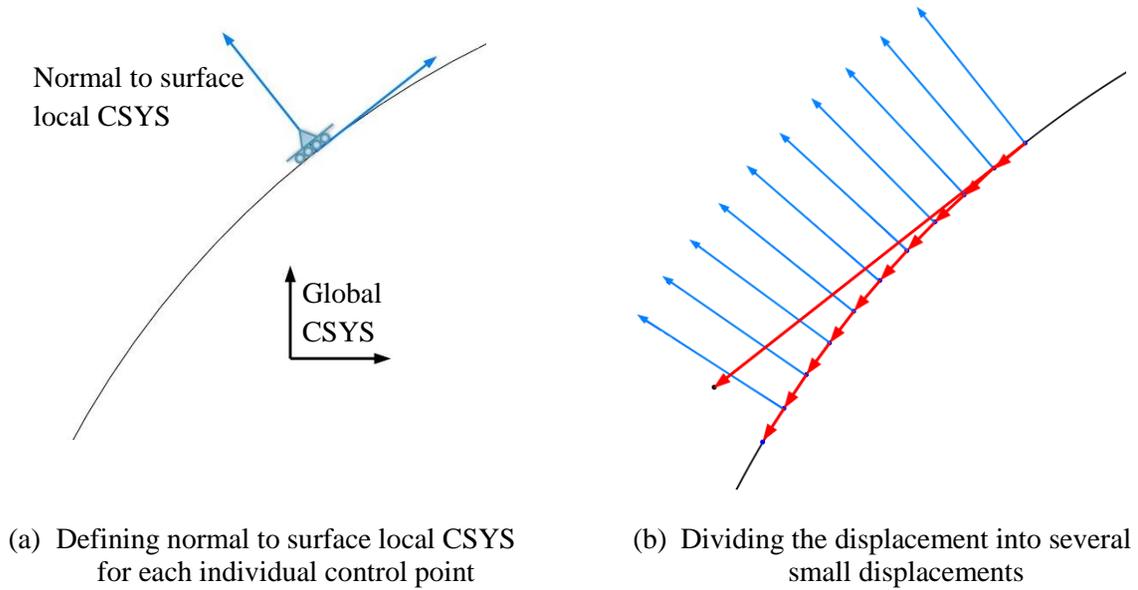

(a) Defining normal to surface local CSYS for each individual control point

(b) Dividing the displacement into several small displacements

Fig. 4. Sliding constraint technique

### 3.4. The solving procedure

According to the descriptions in Section 3, Fig. 5 illustrates the flow diagram of the present multi-step IIGA methodology. In this procedure, the first estimation of the middle configuration is generated using the presented mapping technique. Then, in-plane forces are calculated using the stiffness matrices and the displacements of control points in the local coordinates of $C^F$ and $C^M$. In the subsequent stage, the displacements of the control points are calculated by solving Eq. (7) on the final configuration. Using the continuum relationships described in Ref [47], the principal strains ($\varepsilon_1^{FM}, \varepsilon_2^{FM}, \varepsilon_3^{FM}$) between middle and final forming stages for all elements are calculated. The one-step IIGA is subsequently applied to the updated middle configuration and principal strains between $C^M$ and the initial blank ($\varepsilon_1^{MB}, \varepsilon_2^{MB}, \varepsilon_3^{MB}$) are obtained. Therefore, the total strains in $C^F$ are calculated as below:

$$\varepsilon_i^F = \varepsilon_i^{FM} + \varepsilon_i^{MB} \; ; \; i = 1,2,3 \tag{21}$$

It should be noted that in the first iteration of this method, the thickness and the plastic properties of all elements are assumed to be the thickness of the initial sheet and the elastic properties, respectively. In the next iterations, by calculating equivalent strain and stress for each element, the plastic property matrix is modified as the following equation [48]:

$$[D_p] = \left(\frac{1+r}{1+2r}\right)\frac{\bar{\sigma}}{\bar{\varepsilon}} \begin{bmatrix} 1+r & r & 0 \\ r & 1+r & 0 \\ 0 & 0 & 0.5 \end{bmatrix} \tag{22}$$

In which, r, $\bar{\sigma}$ and $\bar{\varepsilon}$ are the Lankford coefficient, equivalent stress, and equivalent strain, respectively.

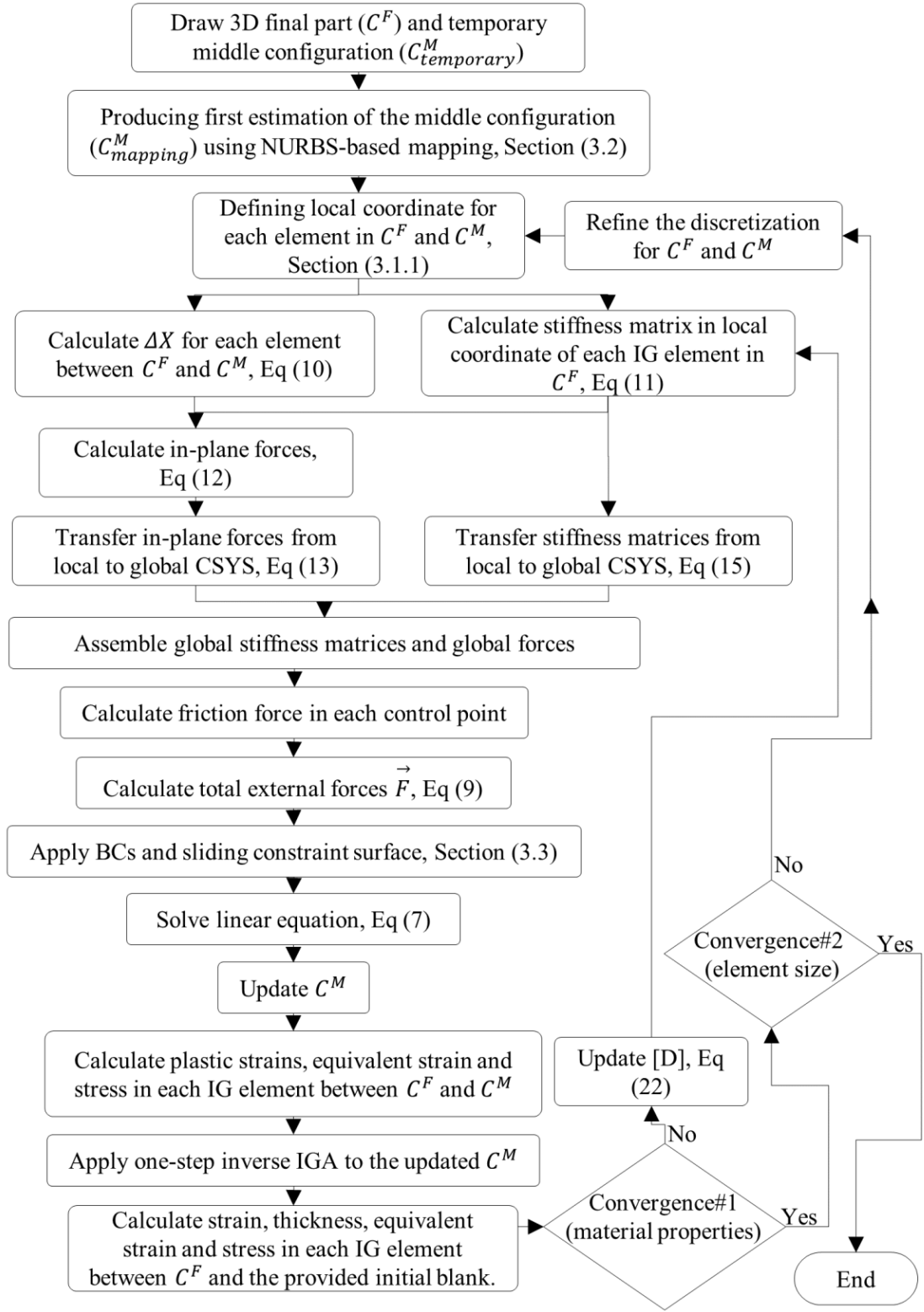

Fig. 5. Flow diagram of the developed multi-step IIGA

As shown in Fig. 5, the present multi-step IIGA solver consists of two different convergence loops. In the first loop, material properties of elements are updated to obtain the approximate results. Then, the number of elements is increased to calculate more accurate results.

## 4. Results and evaluations

To evaluate the credibility of the presented multi-step methodology, two classical problems have been investigated. In the first problem, deep drawing of a rectangular box has been experimentally carried out, and in the second problem two-step drawing of a circular workpiece has been conducted using forward FEM. The results of these processes have been compared to those of the present multi-step IIGA method. In these problems, two different objectives have been taken into consideration. The first objective followed in problem 1 is improving the accuracy of the one-step inverse method by employing the deformation theory of plasticity in more than one stage, and the second objective followed in problem 2 is considering tool changes during a multi-stage forming process. To speed up the computations of the IIGA methodology, the present study uses the second-order NURBS basis functions for solving the problems.

### 4.1. Forming of a rectangular box

Fig. 6 shows the tool geometry parameters used in the experimental setup. An initial square blank ($128mm \times 128mm$) has been formed with a depth of $30\ mm$ by employing a 30-ton hydraulic press at the velocity of 5mm/s. The material properties of the studied sheet and forming parameters are given in Table 1.

Table 1. Material properties of the used sheet and forming parameters of the experimental test

| Parameter | Value |
|---|---|
| Strength coefficient ($K$) | 545 MPa |
| Strain hardening exponent ($n$) | 0.2562 |
| Anisotropy ($r_m$) | 1.1 |
| Young's modulus ($E$) | 206 GPa |
| Poisson ratio ($v$) | 0.3 |
| Friction coefficient ($\mu$) | 0.1 |
| Blank holder force | 10 KN |

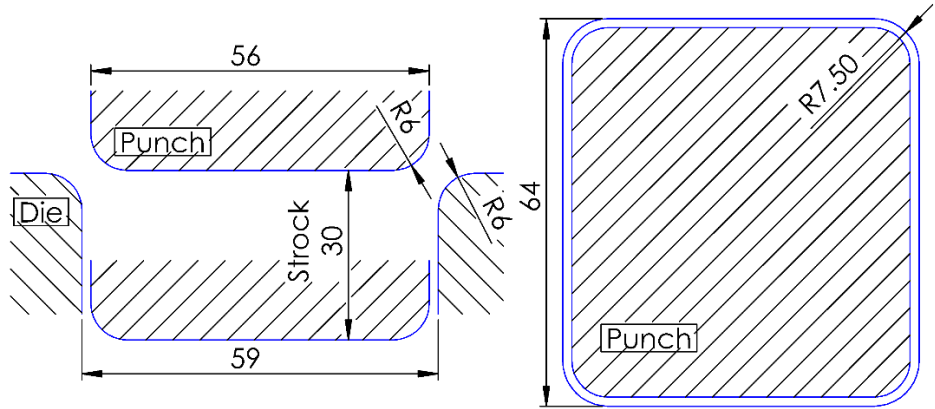

Fig. 6. Dimensions of die and punch in the experiment (All dimensions in mm)

After conducting the experimental test, the final part has been digitized using a coordinate measuring machine (CMM). Also, the thicknesses have been measured along a specific path using a calibrated micrometer. The final part acquired from the experimental setup is shown in Fig. 7 (a). Since the rectangular box has two planes of symmetry, only a quarter of the final part shown in Fig. 7 (b) has been analyzed. This part has been considered as an input of the one-step and multi-step IIGA methods. In IIGA models, the final part has been discretized by $38 \times 38$ control points in the convergence state.

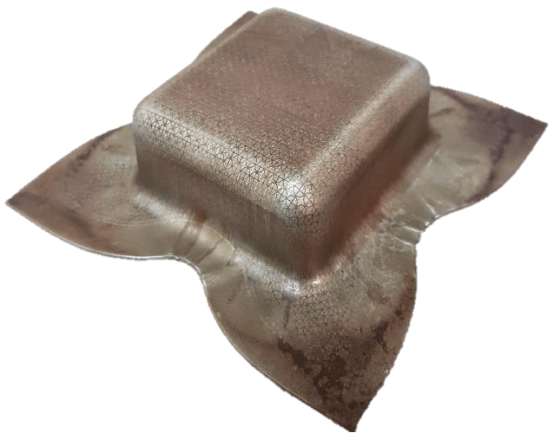 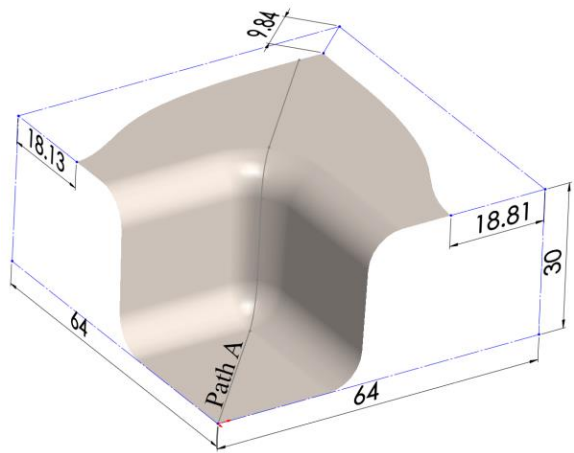

(a) The final part acquired from the experiment  (b) A quarter of the final part

Fig. 7. The studied deep-drawn part in problem 1

Fig. 8 illustrates the user-defined middle part ($C_{user}^M$), the middle part resulted from the NURBS-based mapping ($C_{mapping}^M$), and the last version of the middle part achieved after convergence of the multi-step forming solver. Fig. 9 compares the initial blanks resulted from the one-step IIGA, multi-step IIGA, and the experiment. The maximum errors achieved in the initial blank prediction using one-step and multi-step IIGA models are approximately 2% and 1.4%. It means the error of initial blank estimation in the multi-step IIGA model decreases, compared to that in the one-step IIGA model.

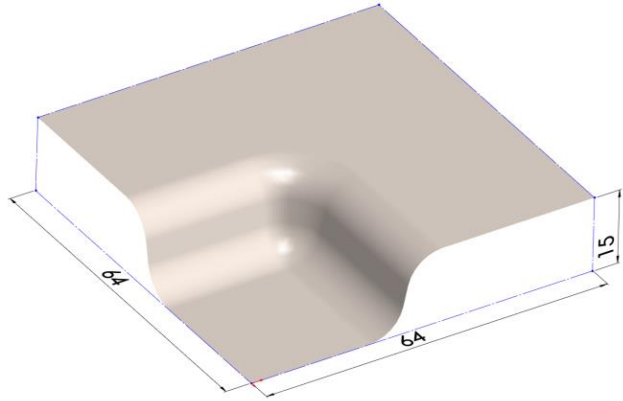

(a) The user-defined middle configuration ($C_{user}^M$)

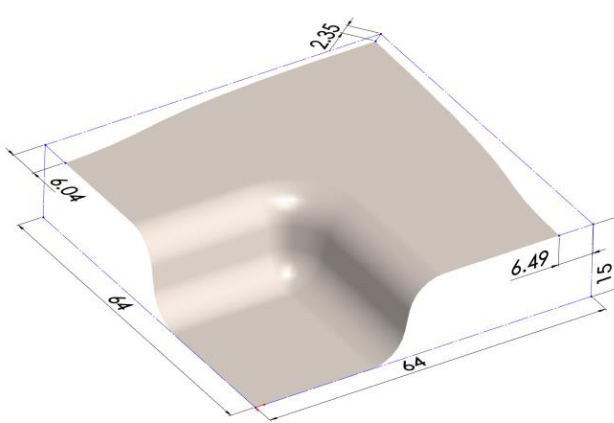 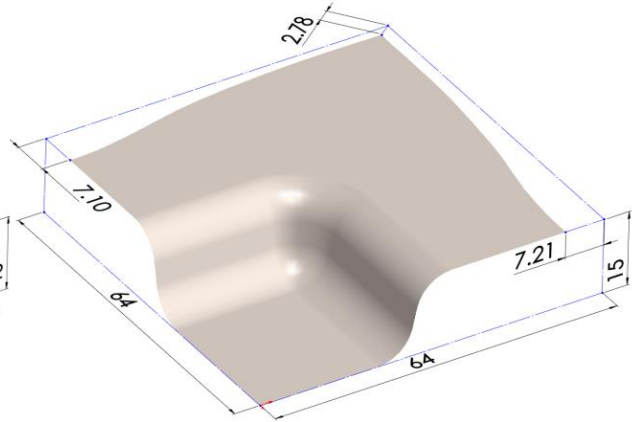

(b) The middle part resulted from the NURBS-based mapping ($C_{mapping}^M$)

(c) The last version of the middle part achieved after the convergence ($C^M$)

Fig. 8. The middle part geometries in problem 1 (All dimensions in mm)

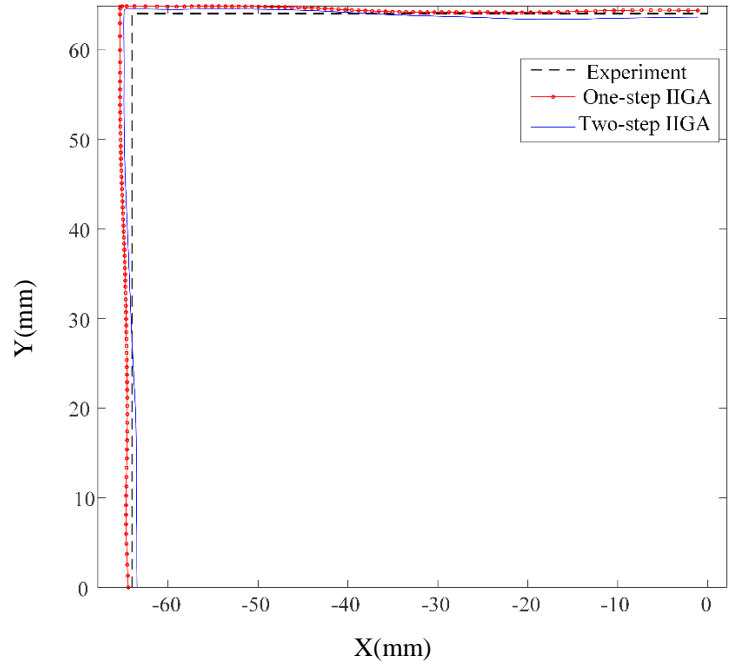

Fig. 9. The calculated initial blanks by one-step IIGA, multi-step IIGA, and experiment

Fig. 10 displays the principal strains calculated by one-step and multi-step IIGA models on the forming limit diagram (FLD), which is created based on the relationships presented in Ref. [49]. Although the multi-step IIGA model predicts no risk of necking, the one-step approach shows necking and fracture during this process. The results of the experiment indicate that necking is not generated during this process. Therefore, the one-step IIGA model cannot accurately predict the forming severity.

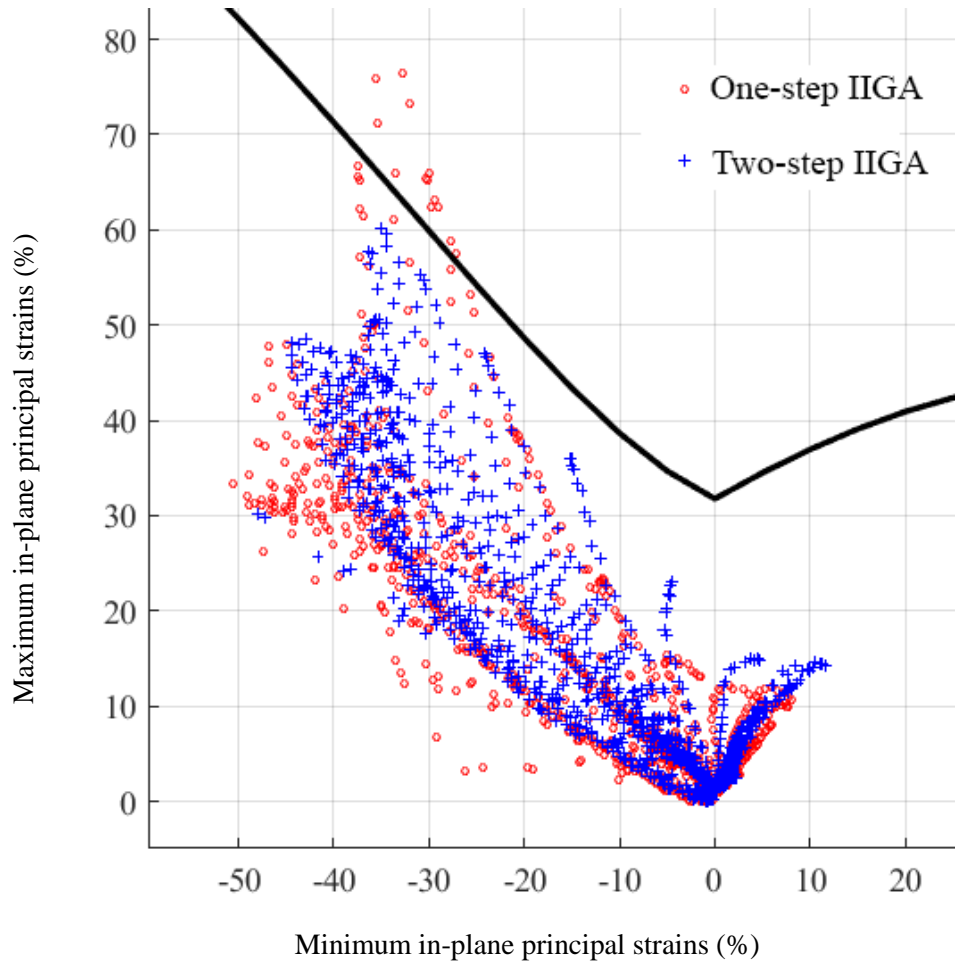

Fig. 10. Illustration of the forming severity of the rectangular box in problem 1

Fig. 11 shows the thickness strain distributions calculated by experiment, forward FEM, one-step IIGA, and multi-step IIGA models, along path A, which is shown in Fig. 7(b).

As revealed in Figs. 9-11, the multi-step IIGA model is able to precisely analyze this severe process, in contrast to the one-step IIGA model, which is not capable of accurately predicting the final part's formability. The main reason of this ability in the multi-step model can be related to employing the deformation theory of plasticity in the initial-middle and middle-final forming operations, in contrast to the one-step model which uses the deformation theory of plasticity only in the initial-final forming stage.

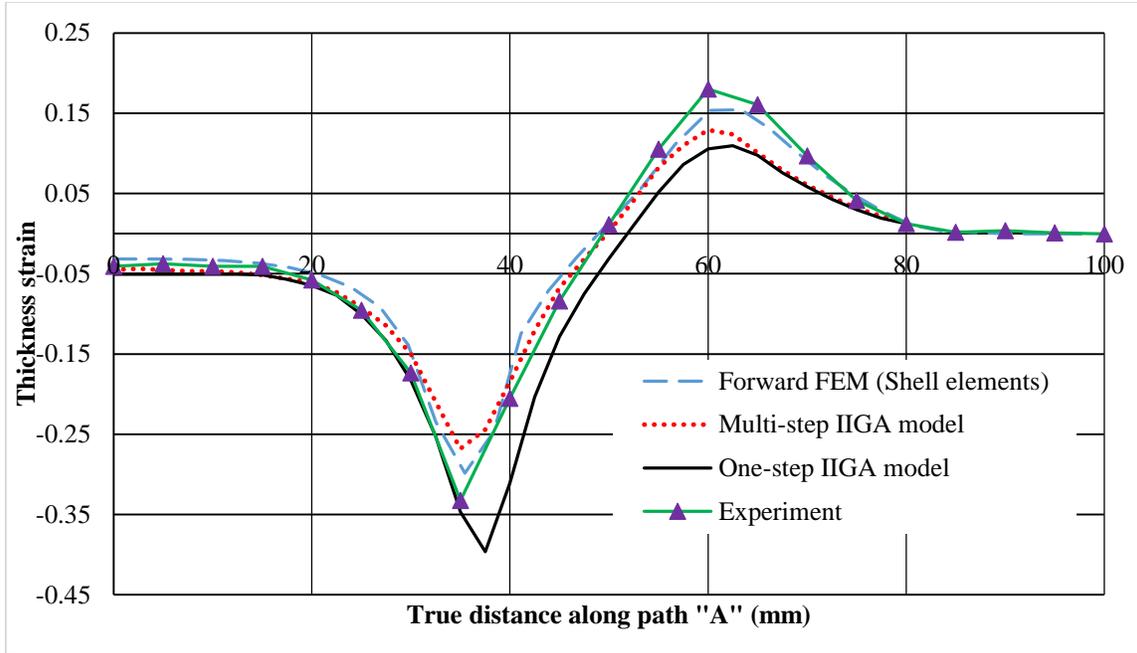

Fig. 11. Comparison of thickness strains along path "A"

## 4.2. Two-step drawing of a circular cup

As the second example, the two-step forming process of a circular part has been studied. The geometry parameters of the first and second punch, die and blank holder are presented in Fig. 12. In this two-step forming process, the final part shown in Fig. 13 (a) has been constructed by two punch movements; punch 1 and then punch 2. This part has been imported into the one-step and multi-step IIGA solvers. The material properties of the used sheet and forming parameters are given in Table 1. The user-defined middle part ($C_{user}^M$), the middle part generated by the NURBS-based mapping ($C_{mapping}^M$), and the middle part obtained in the final step of the multi-step solver are illustrated in Fig. 13 (b, c, d).

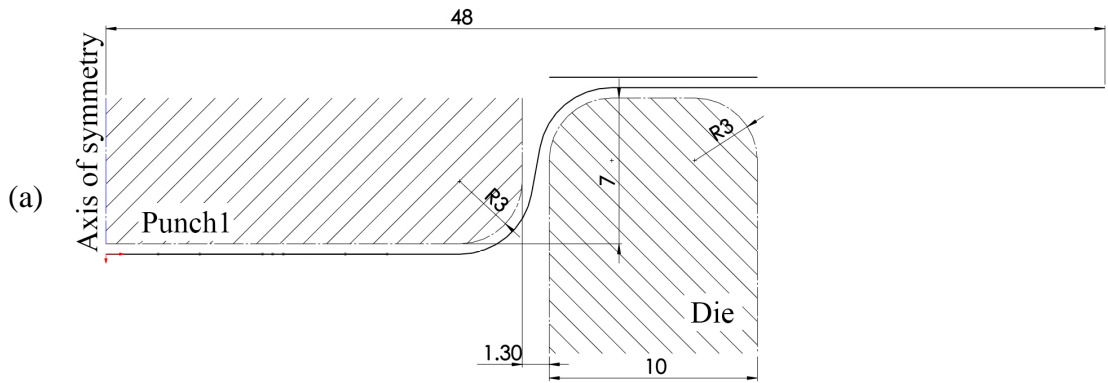

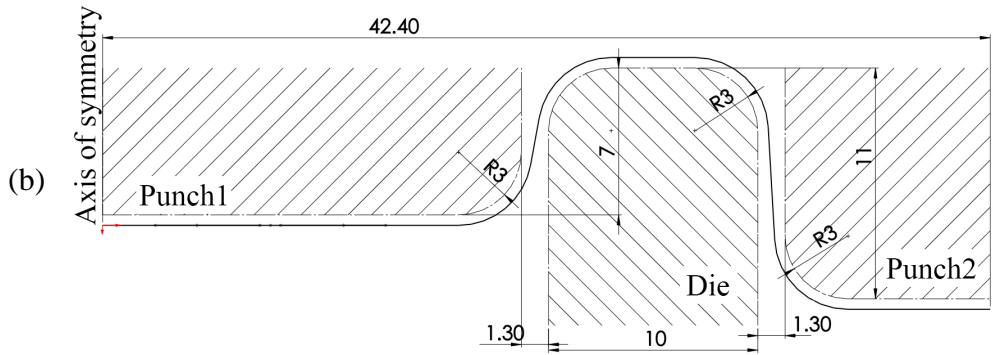

Fig. 12. Forming parameters in problem 2: (a) the first stage (b) the second stage (All dimensions in mm).

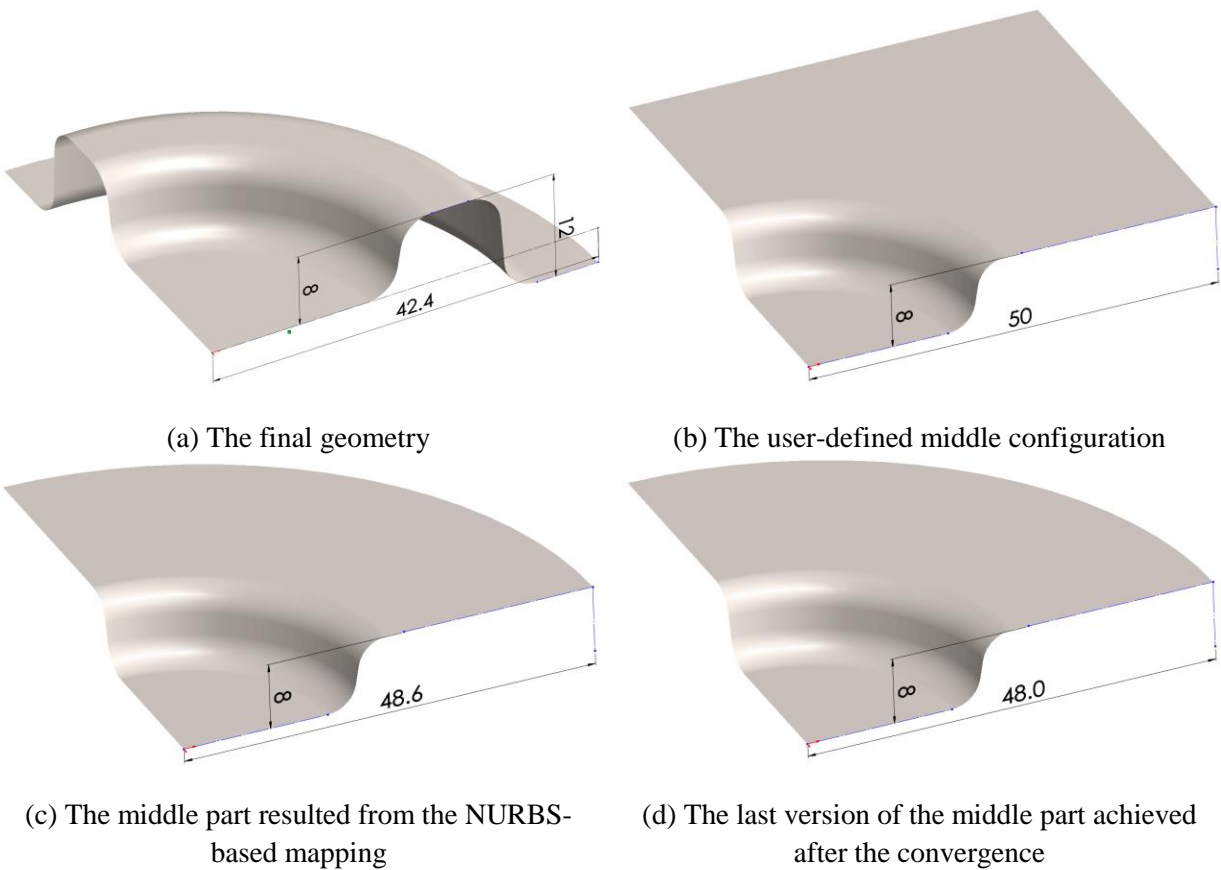

(a) The final geometry

(b) The user-defined middle configuration

(c) The middle part resulted from the NURBS-based mapping

(d) The last version of the middle part achieved after the convergence

Fig. 13. The middle part geometries in problem 2

Thickness strains along the radial path of the final and middle configurations for one-step IIGA, multi-step IIGA, and forward FEM are presented in Fig. 14. As indicated in this figure, the major errors of the one-step and multi-step IIGA models are observed in the zones of punch and die profile radii. This is mainly because the used IIGA models employ membrane elements and cannot predict the bending effects, which are significantly crucial in these regions. Furthermore, the multi-step approach improves the accuracy of the one-step model in the punch and flange zones. These enhancements are achieved by including the second stage of the forming process into the computation.

The results of the initial blank diameters calculated by all approaches are compared in Table 2. As shown in Fig. 14, this forming process is not severe as compared with problem 1. Therefore, a slight difference between the initial blanks predicted by the one-step and multi-step IIGA models is observed.

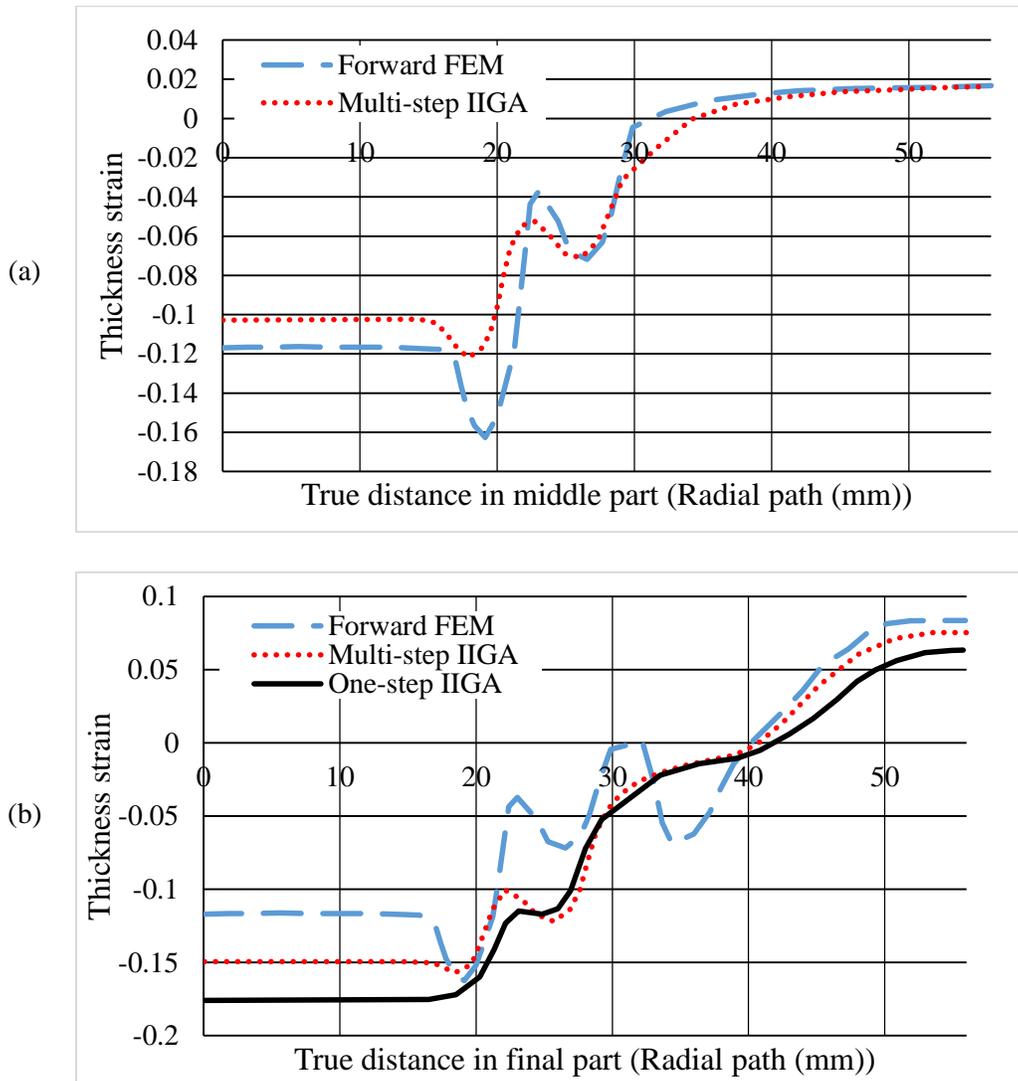

Fig. 14. Radial thickness strain: (a) middle configuration and (b) final part

Table 2: Comparison of the calculated initial blanks

| Parameter | One-step IIGA | Two-step IIGA | Forward FEM |
|---|---|---|---|
| Initial blank diameter (mm) | 99.19 | 99.31 | 100 |
| Error (%) compared with FEM | 0.81 | 0.69 | - |

Although the CPU time required for solving this problem using the presented multi-step IIGA model is about 6 times greater than the time needed in one-step method, the multi-step model requires much less computation time than forward FEM. On the other hand, the multi-step IIGA approach improves the accuracy of results obtained by one-step method. As a result, the best compromise between accuracy and computation time can be achieved by the presented multi-step IIGA solver for analysis of this forming process.

## 5. Conclusion

In this paper, the transfer-based multi-step IIGA model has been presented for more accurate prediction of the initial blank and strains in sheet metal stamping processes. The capability of the presented model has been experimentally and numerically demonstrated in forming of a rectangular box and also two-stage drawing of a circular workpiece. In the present multi-step methodology, a new NURBS-based mapping and sliding constraint technique have been developed to obtain the governing equations which are solved with no concern about the convergence. Since this model uses NURBS basis functions for both exact geometric modeling and inverse forming analysis, the time-consuming process of mesh generation has been eliminated. This important issue leads to fast prediction of forming parameters at the initial stages of design. Therefore, the present model can be successfully used in simulation of sheet metal forming processes.